\documentclass[12pt]{article}
\usepackage[cp1251]{inputenc}
\usepackage{amssymb,latexsym,amsmath}

\usepackage{amsfonts}
\usepackage[english]{babel}
\usepackage{indentfirst}
\usepackage{graphicx,graphics,verbatim}
	\usepackage{cite}

\begin{document}

		
		MSC 68-04

\begin{center}		
	{\bf The explicit formula for solution of anomalous
		diffusion equation in the multi-dimensional space} \\
\end{center}	

\begin{center}
		\textbf{Durdimurod K. Durdiev}\\
	Bukhara Branch of the Institute of Mathematics at the Academy of Sciences of the Republic of Uzbekistan, Bukhara, Uzbekistan\\
	durdiev65@mail.ru\\
		\textbf{Elina L. Shishkina}\\
	Voronezh State University, Voronezh, Russia\\
	shishkina@amm.vsu.ru\\
	\textbf{Sergei M. Sitnik}\\
	Belgorod State University, Belgorod, Russia\\
	sitnik@bsu.edu.ru\\

\end{center}






{\bf Abstract}. This paper intends on obtaining the explicit solution
of $n$-dimensional anomalous diffusion equation in the infinite
domain  with non-zero  initial condition and vanishing condition
at infinity. It is shown that this equation can be derived from
the parabolic integro-differential equation with memory in which
the kernel is $t^{-\alpha}E_{1-\alpha, 1-\alpha}(-t^{1-\alpha}), \
\alpha\in(0, 1),$ where $E_{\alpha, \beta}$ is the Mittag-Liffler
function. Based on Laplace and Fourier transforms the properties
of the Fox H-function and convolution theorem, explicit solution
for anomalous diffusion equation is obtained.


\tableofcontents

\section{Introduction to the problem and its setting}

The processes of transport and aftereffect in disordered inhomogeneous composite media often have kinetics
that do not obey the normal (Gaussian) distribution. For this reason, such processes are usually called anomalous  or non-classical.
They are observed experimentally in the study of diffusion in turbulent flows, heat and mass transfer in plasma, filtration of fluids
in inhomogeneous porous media, impurity transfer in complex geological formations, charge transfer in amorphous semiconductors,
relaxation processes in polymers, hydrodynamics of non-Newtonian fluids, the evolution of complex biological systems, the transfer
of information resources to global communication networks and many other phenomena.

The rapid development of fractional differential equations was largely due to the
discovered practical applications of fractional calculus, primarily in the physics
of complex inhomogeneous media. It turned out that fractional differential equations
are ideally suited for modeling anomalous processes occurring in systems with a fractal
structure or having a power-law memory.  As a result, integro-differentiation of fractional order began to develop
as a powerful modern apparatus of mathematical modeling, including both analytical and numerical methods for studying
fractional differential models.

The most common and studied fractional differential models are models of various anomalous diffusion-type transfer processes,
the kinetics of which are well described (at least asymptotically) by power laws with fractional exponents.
The value of this exponent, depending on the standard deviation of the positions of the diffusing particles
(if its final value exists) versus time, allows us to divide the processes of anomalous diffusion into subdiffusion ones
and superdiffusion ones.  With subdiffusion, the particle transfer
intensity is lower than with classical diffusion. Such processes are usually observed in systems with memory. The superdiffusion
phenomenon has a higher particle transfer rate than classical diffusion and is often observed in media with a fractal structure.
In systems with both spatial and temporal non locality, both sub-  and  superdiffusion  processes can occur \cite{1}.
Examples of fractional differential equations of the diffusion type are the properly equations of subdiffusion and
superdiffusion \cite{2,3,4,5,6}, diffusion-wave equations \cite{7, 8, 9, 10, 11}, fractional
differential equations of Fokker-Plann \cite{12, 13, 14, 15, 16}, advection-dispersion \cite{17, 18, 19, 20},
reaction-diffusion \cite{21, 22, 23} and a number of others.

The solvability of Cauchy problems and initial-boundary value problems for various types
of  linear fractional differential equations of diffusion type were investigated in the
works of A.A. Kilbas, S.M. Sitnik, A.N. Kochubey, Yu. Luchko,  R.  Gorenflo, F.  Mainardi,
M. M. Meerschaert, G. Pagnini, J.J. Trujillo and many others \cite{24, 25, 26, 27, 28, 29, 30, 31, 32, 33, 34}. To construct a solution of
linear fractional differential equations of diffusion type,
various methods and algorithms based on the Green's function,
Fourier, Laplace, and Mellin integral transforms, a generalization
of the method of separation of variables, reduction to
Volterra-type integral equations, and several others were
proposed. At the same time, there are practically no methods for
obtaining analytical solutions of fractional differential
equations describing anomalous diffusion processes.

In this paper, we consider the following $n-$dimensional
integro-differential equation of the anomalously diffusive
transport of solute in heterogeneous porous media \cite{35}
$$u_t(x,t)+\,_0^CD_t^{\alpha}u-\Delta{u(x,t)}=f(x,t), \eqno{(1.1)}$$
which satisfies the initial and boundary conditions
$$u(x,0)=g(x), \qquad \lim_{|x|\rightarrow\infty}\left(u, \nabla u\right)(x,t)=0, \qquad t>0, \qquad x=(x_1,x_2,...,x_n)\in\mathbb{R}^n, \eqno{(1.2)}$$
where the Gerasimov--Caputo fractional differential operator
$_0^CD_t^{\alpha}$ is defined by \cite{25}
$$\,_0^CD_t^{\alpha}g(t):=_0I_t^{1-\alpha}f'(t)=\frac{1}{\Gamma(1-\alpha)}\int\limits_0^t\frac{f'(\tau)}{(t-\tau)^{\alpha}}d\tau,$$
where
$$_0I_t^{\alpha}g(t):=\frac{1}{\Gamma(\alpha)}\int\limits_0^t\frac{f(\tau)}{(t-\tau)^{1-\alpha}}d\tau,$$
$\Delta$ is the $n-$dimensional Laplace operator with respect to $x$ and  $\nabla=\left(\frac{\partial}{\partial x_1},...,\frac{\partial}{\partial x_n}\right).$

\section{Preliminaries}

In this section, we present well known definitions and lemmas that will be used for proof of main results.

The Mittag-Leffler functions $E_{\alpha}(z)$ and $E_{\alpha, \beta}(z)$ are defined by the following series:
$$E_{\alpha}(z):=\sum_{n=0}^{\infty}\frac{z^n}{\Gamma(\alpha{n}+1)}=:E_{\alpha,1}(z)$$
and
$$E_{\alpha, \beta}(z):=\sum_{n=0}^{\infty}\frac{z^n}{\Gamma(\alpha{n}+\beta)},$$
respectively, where $\alpha,z, \rho\in\mathbb{C}; \mathfrak{R}(\alpha)>0.$  These functions are natural extensions of the exponential, hyperbolic and trigonometric functions, since
$$ E_1(z)=e^z,\quad  E_2(z^2)={\rm cosh}\,z,\quad  E_2(-z^2)=\cos z,\quad E_{1,2}(z)=\frac{e^z-1}{z},\quad E_{2,2}(z^2)=\frac{{\rm sinh}\,z}{z}.$$
The three-parameter Mittag-Leffler function or Prabhakar function
is \cite{36}:
$$E_{\alpha, \beta}^{\gamma}:=\sum\limits_{n=0}^{\infty}\frac{(\gamma)_n}{\Gamma(\alpha{n}+\beta)}\frac{z^n}{n!}, \eqno{(2.1)}$$
where $\alpha, \beta,\gamma, z\in\mathbb{C}$, and $(\gamma)_n$ denotes the Pochammer symbol or the shifted factorial defined by
$$(\gamma)_0=1, \qquad  (\gamma)_n=\gamma(\gamma+1)...(\gamma+n-1), \qquad   \gamma\neq{0}.$$
Also we can write
$$(\gamma)_n=\frac{\Gamma(\gamma+n)}{\Gamma(\gamma)}, $$
where $\Gamma(\gamma)$ is the Gamma function.
We have  following special cases:
$ E_{\alpha, \beta}^1(z)=E_{\alpha, \beta}(z)$ and $ E_{\alpha,1}^1=E_{\alpha}(z).$

Recall  that the function (2.1) can be rewritten in terms of the Fox H-function as \cite{36, 37}:
$$E_{\alpha,\beta}^{\gamma}(z)=\frac{1}{\Gamma(\gamma)}H_{1,2}^{1,1}
\left[-z\left|{^{(1-\gamma,1)}_{(0,1),    (1-\beta,\alpha)}}\right.\right]$$
We define the integral operator $\mathrm{\mathcal{E}}_{\alpha, \beta, \omega; a+}^{\gamma}$ as follows \cite{36,38}:
$$\left(\mathrm{\mathcal{E}}_{\alpha,\beta,\omega; a+}^{\gamma}\varphi\right)(t):=\left(t^{\beta-1}E_{\alpha, \beta}^{\gamma}(\omega{t}^{\alpha})\right)*\varphi(t)=\int\limits_a^t(t-\tau)^{\beta-1}E_{\alpha,\beta}^{\gamma}(\omega(t-\tau)^{\alpha})\varphi(\tau)d\tau. \eqno{(2.2)}$$
Note the integral operator (2.2) is nowadays known in literature as Prabhakar fractional integral.

\textbf{Lemma 1.} The following  Laplace transform of a three-parameter Mittag-Leffler function  is true \cite{36,39}:
$$L\left[t^{\beta-1}E_{\alpha,\beta}^{\gamma}(\pm{\omega}t^{\alpha})\right](s)=
\int\limits_0^{\infty}e^{-st}t^{\beta-1}E_{\alpha,\beta}^{\gamma}(\pm{\omega}t^{\alpha})dt=\frac{s^{\alpha\gamma-\beta}}{(s^{\alpha}\mp{\omega})^{\gamma}},$$
where $\left|\omega/{s^{\alpha}}\right|<1.$

\textbf{Lemma 2.} The Laplace transform of $e^{-\lambda{t}}t^{\beta-1}E_{\alpha,\beta}^{\gamma}(\pm\omega{t}^{\alpha})$ is given by the following formula \cite{39}:
$$L\left[e^{-\lambda{t}}t^{\beta-1}E_{\alpha\beta}^{\gamma}(\pm\omega{t}^{\alpha})\right](s)=
\frac{(s+\lambda)^{\alpha\gamma-\beta}}{((s+\lambda)^{\alpha}\mp\omega)^{\gamma}}, $$
where $\lambda\geq{0}, $  $\left|\omega/(s+\lambda)^{\alpha}\right|<1.$

In the case $\lambda=0,$ Lemma 2 coincides with Lemma 1.

\textbf{Lemma 3.} For arbitrary $\alpha>0,$ $ \beta$ is an arbitrary complex number, $\mu>0$ and $a\in\mathbb{R}$, the following formula is valid \cite{20}:
$$\int\limits_{\mathbb{R}^n}e^{i\xi\cdot{x}}E_{\alpha,\beta}^{(n)}(-a|\xi|^{\mu})d\xi=(2\pi)^{n/2}|x|^{1-n/2}\int\limits_0^{\infty}|\xi|^{n/2}E_{\alpha,\beta}^{(n)}(-a|\xi|^{\mu}){J}_{\frac{n}{2}-1}(|x||\xi|)d|\xi|.$$
Here ${J}_{\frac{n}{2}-1}(\cdot)$ is a Bessel function and $ E_{\alpha,\beta}^{(n)}(z)$ denotes
$n-$th derivatives of the Mittag-Leffler function. $n-$th derivatives of the Mittag-Leffler function can be expressed in terms of the Fox H-function as
$$E_{\alpha,\beta}^{(n)}(z)=H_{1,2}^{1,1}\left[-z\left|{ ^{(-n,1)}_{(0,1),   (1-(\alpha{n}+\beta),\alpha)}}\right.\right].$$

\textbf{Lemma 4.} If $\{k(t),r(t)\}\in{L_1[0,T]}$ for a fixed $T>0$
and $k(t), r(t)$ are connected by the integral equation
$$r(t)=k(t)+\int\limits_0^tk(t-\tau)r(\tau)d\tau,\qquad  t\in[0,T], \eqno{(2.3)}$$
then  the solution of the integral equation
$$\varphi(t)=\int\limits_0^tk(t-\tau)\varphi(\tau)d\tau+f(t),\qquad f(t)\in{L_1[0,T]}\eqno{(2.4)}$$
is expressed by formula
$$\varphi(t)=\int\limits_0^tr(t-\tau)f(\tau)d\tau+f(t). \eqno{(2.5)}$$

Although the assertion of Lemma 4 is well known, we give a
original proof.

\textbf{Proof} of Lemma 4.  Let the equality  (2.3) be satisfied.
Then  from (2.3) we have
$$k(t)=r(t)-\int\limits_0^tr(t-\tau)k(\tau)d\tau.$$
Taking into account this relation, from (2) we obtain the
following chain of equalities:
$$\varphi(t)=\int\limits_0^t\left[r(t-\tau)-\int\limits_0^{t-\tau}r(t-\tau-\alpha)k(\alpha)d\alpha\right]\varphi(\tau)d\tau+f(t)=\int\limits_0^tr(t-\tau)\varphi(\tau)d\tau-$$
$$-\int\limits_0^{t}\left(\int\limits_0^{t-\tau}r(t-\tau-\alpha)k(\alpha)d\alpha\right)\varphi(\tau)d\tau+f(t)=
\int\limits_0^tr(\tau)\varphi(t-\tau)d\tau-$$
$$-\int\limits_0^{t}\left(\int\limits_0^{t-\tau}r(\alpha)k(t-\tau-\alpha)d\alpha\right)\varphi(\tau)d\tau+f(t)=\int\limits_0^tr(\tau)\varphi(t-\tau)d\tau-$$
$$-\int\limits_0^{t}\left(\int\limits_0^{\tau}r(\alpha)k(\tau-\alpha)d\alpha\right)\varphi(t-\tau)d\tau+f(t)=\int\limits_0^tr(\tau)\varphi(t-\tau)d\tau-$$
$$-\int\limits_0^{t}r(\alpha)\left(\int\limits_{\alpha}^{t}k(\tau-\alpha)\varphi(t-\tau)d\tau\right)d\alpha+f(t)=$$
$$=\int\limits_0^tr(\tau)\left[\varphi(t-\tau)-\int\limits_0^{t-\tau}k(\tau)\varphi(t-\tau-\alpha)d\alpha\right]d\tau+f(t).$$

According to (2.4), the expression  in the square brackets of the
right side in the last equalities is equal to $f(t)$. Therefore,
we get (2.5).

Now we  show the assertion of Lemma 4 is true  in the opposite
direction. Indeed, from (2.5) on based of (2.3), we obtain

$$\varphi(t)=\int\limits_0^t\left[k(t-\tau)+\int\limits_0^{t-\tau}k(t-\tau-\alpha)r(\alpha)d\alpha\right]f(\tau)d\tau+f(t)=$$
$$=\int\limits_0^tk(t-\tau)f(\tau)d\tau+\int\limits_0^t\left(\int\limits_0^{t-\tau}k(t-\tau-\alpha)r(\alpha)d\alpha\right)f(\tau)d\tau+f(t)=$$
$$=\int\limits_0^tk(\tau)\left[f(t-\tau)+\int\limits_0^{t-\tau}r(\alpha)f(t-\tau-\alpha)d\alpha\right]d\tau+f(t).$$
From this, in view of (2.5), we obtain (2.3).

\section{Explicit solution of the problem (1.1) and (1.2)}

Let
$$\widetilde{f}(\xi,t):=\int\limits_{\mathbb{R}^n}f(x,t)e^{-i\xi\cdot{x}}dx,$$
$$
x\in\mathbb{R}^{n},\qquad \xi\in\mathbb{R}^{n},\qquad \xi\cdot{x}=\sum_{j=1}^{n}\xi_{i}\cdot{x_i}, \qquad dx=dx_1dx_2...dx_n
$$ is the Fourier transform of $f(x,t)$ with respect to the spatial variable $x$ and
$$
\widetilde{g}(\xi):=\int\limits_{\mathbb{R}^n}g(x)e^{-i\xi\cdot{x}}dx.$$
The Laplace transform of a function $u(x,t)$ with respect to the variable $t$ is given by
$$
\hat{u}(x,s)=\int\limits_0^{\infty}e^{-st}u(x,t)dt,\qquad t>0.
$$

The unknown function $u(x,t)$ is required to be sufficiently well behaved to be treated with its derivatives $u_t(x,t),$ $u_{x_ix_i}(x,t), i=1,\ldots, n$ by technique of Laplace (in $t$)
and Fourier (in $x$) transforms. The given  functions $f(x, t)$ and $g(x)$ are also assumed to have such properties.

\textbf{Theorem 1.} The explicit solution of the problem (1.1) and (1.2) can be expressed by formula

$$u(x, t)=\frac{1}{(2\pi)^{n}}\int\limits_{\mathbb{R}^n}\sum_{j=0}^{\infty}(-1)^j\left(\mathcal{E}_{1, (1-\alpha)j+1,-|\xi|^2; 0+}^{j+1}\tilde{f}\right)(\xi,\tau)e^{i\xi\cdot x}d \xi+$$$$+\int\limits_{\mathbb{R}^n}G(x-\xi, t)g(\xi)d \xi, \qquad d\xi=d\xi_1d\xi_2...d\xi_n, \eqno(3.1)$$
where the Green function $G(x, t)$ is given by
$$G(x, t)=\frac{1}{(2\pi)^{n/2}|x|^n}\sum_{j=0}^{\infty}\frac{\left(-t^{1-\alpha}\right)^j}{j!}\times$$
$$\times
\left\{H_{2,0}^{1,2}\left[\frac{|x|^2}{2t^{1/2}}\left|
\begin{array}{c}
$$(1+(1-\alpha)j, 1/2)$$ \\
$$(n/2, 1/2), (2+j, 1/2) $$ \\
\end{array}\right.\right]+H_{2,0}^{1,2}\left[\frac{|x|^2}{2t^{1/2}}\left|
\begin{array}{c}
$$(1+(1-\alpha)j, 1/2)$$ \\
$$(n/2, 1/2), (1+j, 1/2) $$ \\
\end{array}\right.\right]\right\},
$$
where
$x\in\mathbb{R}^n,$ $\xi\in\mathbb{R}^n.
$

\textbf{Proof.}
Let $F\{u(x,t)\}:=\tilde{u}(\xi,t)$ be the Fourier transform of $u(x,t)$ with respect to variable $x$, and $L\{u(x,t)\}:=\hat{u}(x,s)$ be the Laplace transform of $u(x,t)$ with respect to variable $t$. In sequence, applying to the equation (1.1) the Laplace transform with respect to the time variable $t$  and the Fourier transform with respect to the spatial variable $x$, we obtain the following equation:
$$ s\hat{\tilde{u}}(\xi,s)-\tilde{u}(\xi,0)+s^{\alpha-1}\left[s\hat{\tilde{u}}-\tilde{u}(\xi,0)\right]=-|\xi|^2\hat{\tilde{u}}(\xi,s)+\hat{\tilde{f}}(\xi,s), \xi\in\mathbb{R}^n.  $$

Taking into account the initial condition (1.2), this equation  yields
$$\hat{\tilde{u}}(\xi,s)=\frac{1}{s+s^{\alpha}+|\xi|^2}\hat{\tilde{f}}(\xi,s)+\frac{1+s^{\alpha-1}}{s+s^{\alpha}+|\xi|^2}\tilde{g}(\xi). \eqno{(3.2)}$$

First, we calculate the inverse Laplace and Fourier transforms of the first term on the  right side of (3.2). It may be performed by using the equality
$$ \frac{1}{s+s^{\alpha}+|\xi|^2}=\frac{s^{-\alpha}}{s^{1-\alpha}+1}\cdot\frac{1}{1+\frac{|\xi|^2s^{-\alpha}}{s^{1-\alpha}+1}} \eqno{(3.3)}$$
and expanding the third factor on the right of this equation into an infinitely decreasing geometric series:
$$\frac{1}{1+\frac{|\xi|^2s^{-\alpha}}{s^{1-\alpha}+1}}=\sum_{n=0}^{\infty}(-|\xi|^2)^n\frac{s^{-\alpha{n}}}{(s^{1-\alpha}+1)^n}$$ for $\left|\frac{|\xi|^2s^{-\alpha}}{s^{1-\alpha}+1}\right|<1.$

In view of (3.3) from last equality we have
$$\frac{1}{s+s^{\alpha}+|\xi|^2}=\sum_{n=0}^{\infty}(-|\xi|^2)^n\frac{s^{-\alpha{(n+1)}}}{(s^{1-\alpha}+1)^{n+1}}.$$

Then, according to Lemma 1, the first term of equation  (3.2) can be expressed as
$$\frac{1}{s+s^{\alpha-1}+|\xi|^2}\hat{\tilde{f}}(\xi,s)=
L\left[\sum_{n=0}^{\infty}(-|\xi|^2)^nt^nE_{1-\alpha,n+1}^{n+1}(-t^{1-\alpha})\right]L\left[\tilde{f}(\xi,t)\right]. \eqno{(3.4)}$$

We now transform the second term on the right side of (3.2). For this we note

$$\frac{1+s^{\alpha-1}}{s+s^{\alpha-1}+|\xi|^2}=(1+s^{\alpha-1})\cdot\frac{s^{-\alpha}}{s^{1-\alpha}+1}\cdot
\frac{1}{1+\frac{|\xi|^2s^{-\alpha}}{s^{1-\alpha}+1}}=$$
$$=\frac{s^{-\alpha}+s^{-1}}{s^{1-\alpha}+1}\sum_{n=0}^{\infty}(-|\xi|^2)^n\frac{s^{-\alpha{n}}}{(s^{1-\alpha}+1)^n}=$$
$$=\sum_{n=0}^{\infty}(-|\xi|^2)^n\frac{s^{-\alpha{(n+1)}}}{(s^{1-\alpha}+1)^{n+1}}+\sum_{n=0}^{\infty}
(-|\xi|^2)^n\frac{s^{-\alpha{n}-1}}{(s^{1-\alpha}+1)^{n+1}}.
$$
In view of last relations, applying Lemma 1 to the second term of (3.2),  we obtain
$$\frac{1+s^{\alpha-1}}{s+s^{\alpha-1}+|\xi|^2}\tilde{g}(\xi)=
L\left[\sum_0^{\infty}(-|\xi|^2)^nt^n(E_{1-\alpha,n+1}^{n}(-t^{1-\alpha})+E_{1-\alpha,n+1}^{n+1}(-t^{1-\alpha})\right]\cdot\tilde{g}(\xi). \eqno{(3.5)}$$

Further, in accordance with the Mittag-Leffler function definition (2.1), from equation (3.4) we get

$$\sum_{n=0}^{\infty}(-|\xi|^2)^nt^n(E_{1-\alpha,n+1}^{n}(-t^{1-\alpha})
=\sum_{n=0}^{\infty}(-|\xi|^2)^nt^n\sum_{j=0}^{\infty}\frac{(n+1)_j}{\Gamma((1-\alpha)j+n+1)}\cdot\frac{(-t^{1-\alpha})^j}{j!}=$$
$$=\sum_{n=0}^{\infty}\sum_{j=0}^{\infty}(-t^{1-\alpha})^j\frac{(j+1)_n}{\Gamma((1-\alpha)j+n+1)}\frac{(-|\xi|^2t)^n}{n!}=
\sum_{j=0}^{\infty}(-t^{1-\alpha})^jE_{1,(1-\alpha)j+1}^{j+1}(-|\xi|^2t).$$

By virtue of this fact we continue converting of the right side of (3.4) as
$$L\left[\sum_{n=0}^{\infty}(-|\xi|^2)^nt^nE_{1-\alpha,n+1}^{n+1}(-t^{1-\alpha})\right]L\left[\tilde{f}(\xi,t)\right]=
L\left[\sum_{j=0}^{\infty}\left(-t^{1-\alpha}\right)^jE_{1, (1-\alpha)j+1}^{j+1}(-|\xi|^2t)\right]\times$$$$\times L\left[\tilde{f}(\xi,t)\right]=
L\left[\left(\sum_{j=0}^{\infty}\left(-t^{1-\alpha}\right)^jE_{1, (1-\alpha)j+1}^{j+1}(-|\xi|^2t)\right)\ast\tilde{f}(\xi,t)\right].$$

Taking into consideration the convolution property of the Laplace transform and the definition of integral operator $\mathrm{\mathcal{E}}_{\alpha,\beta,\omega;a+}^{\gamma}\varphi$ by (2.2), the inverse transform of the first term in equation (3.2) can be obtained as follows:
$$L^{-1}\left[\frac{1}{s+s^{\alpha-1}+|\xi|^2}L\left[\tilde{f}(\xi,t)\right](s)\right]=\sum_{j=0}^{\infty}(-1)^j\left(\mathcal{E}_{1, (1-\alpha)j+1,-|\xi|^2; 0+}^{j+1}\tilde{f}\right)(\xi,\tau). \eqno{(3.6)}$$

Analogically, the inverse Laplace transform of the second term in equation (3.2) by virtue of (3.5) can be expressed as
$$L^{-1}\left[\frac{1+s^{\alpha-1}}{s+s^{\alpha-1}+|\xi|^2}\tilde{g}(\xi)\right]=
\sum_{n=0}^{\infty}\left(-|\xi|^2\right)^nt^n\left(E_{1-\alpha,n+1}^n(-t^{1-\alpha})+E_{1-\alpha,n+1}^{n+1}(-t^{1-\alpha})\right)\tilde{g}(\xi)=$$
$$=\sum_{j=0}^{\infty}(-t^{1-\alpha})^j\left(E_{1,(1-\alpha)j+1}^j(-|\xi|^2t)+E_{1,(1-\alpha)j+1}^{j+1}(-|\xi|^2t)\right)\tilde{g}(\xi).$$

Considering the relationship between the generalized Mittag-Leffler function and the Fox-H function,
the last equality can be rewritten in the form \cite{37}:
$$L^{-1}\left[\frac{1+s^{\alpha-1}}{s+s^{\alpha-1}+|\xi|^2}\tilde{g}(\xi)\right]=\sum_{j=0}^{\infty}(-1)^jt^{(1-\alpha)^j}\times$$$$\times
\left(\frac{1}{\Gamma(j)}H_{1,2}^{1,1}\left[|\xi|^2t\left|^{(1-j,1)}_{(0,1), (-(1-\alpha)j,1)}\right.\right]+\frac{1}{\Gamma(j+1)}H_{1,2}^{1,1}\left[|\xi|^2t\left|^{(-j,1)}_{(0,1), (-(1-\alpha)j,1)}\right.\right]\right)\tilde{g}(\xi)=$$
$$=\sum_{j=0}^{\infty}\frac{(-1)^j}{j!}t^{(1-\alpha)^j}\left(jH_{1,2}^{1,1}\left[|\xi|^2t\left|^{(-j,1)}_{(0,1), (-(1-\alpha)j,1)}\right.\right]+H_{1,2}^{1,1}\left[|\xi|^2t\left|^{(-j,1)}_{(0,1), (-(1-\alpha)j,1)}\right.\right]\right)\tilde{g}(\xi).\eqno{(3.7)}$$

We introduce the following notations:
$$\tilde{h}_{0j}(\xi,{t})=H_{1,2}^{1,1}\left[|\xi|^2t\left|^{(1-j,1)}_{(0,1), (-(1-\alpha)j,1)}\right.\right], \ \ \tilde{h}_{1j}(\xi,{t})=H_{1,2}^{1,1}\left[|\xi|^2t\left|^{(-j,1)}_{(0,1), (-(1-\alpha)j,1)}\right.\right] \eqno{(3.8)}$$
and
$$\tilde{G}(\xi,t)=\sum_{j=0}^{\infty}\frac{(-t^{1-\alpha})^j}{j!}\left[j\tilde{h}_{0j}(\xi,{t})+\tilde{h}_{1j}(\xi,{t})\right]. \eqno{(3.9)}$$

Applying the inverse transform $F^{-1}$ to equations (3.8), we obtain
$$h_{0j}(x,t)=\frac{1}{(2\pi)^n}\int\limits_{\mathbb{R}^n}H_{1,2}^{1,1}\left[|\xi|^2t\left|^{(1-j,1)}_{(0,1), (-(1-\alpha)j,1)}\right.\right]e^{i\xi\cdot x}d\xi, \eqno{(3.10)}$$
$$h_{1j}(x,t)=\frac{1}{(2\pi)^n}\int\limits_{\mathbb{R}^n}H_{1,2}^{1,1}\left[|\xi|^2t\left|^{(-j,1)}_{(0,1), (-(1-\alpha)j,1)}\right.\right]e^{i\xi\cdot x}d\xi. \eqno{(3.11)}$$
Using Lemma 3, we obtain the following results from equations (3.10) and (3.11)
$$h_{0j}(x,t)=\frac{1}{(2\pi)^{n/2}}|x|^{1-\frac{n}{2}}\int\limits_0^{\infty}|\xi|^{n/2}H_{1,2}^{1,1}\left[|\xi|^2t\left|^{(1-j,1)}_{(0,1), (-(1-\alpha)j,1)}\right.\right]\mathfrak{J}_{\frac{n}{2}-1}(|x||\xi|)d|\xi|,$$
$$h_{1j}(x,t)=\frac{1}{(2\pi)^{n/2}}|x|^{1-\frac{n}{2}}\int\limits_0^{\infty}|\xi|^{n/2}H_{1,2}^{1,1}\left[|\xi|^2t\left|^{(-j,1)}_{(0,1), (-(1-\alpha)j,1)}\right.\right]\mathfrak{J}_{\frac{n}{2}-1}(|x||\xi|)d|\xi|.$$

Taking into account a Hankel transform and the properties Fox-H function \cite{37,41}, last two equations can be written as
$$h_{0j}(x,t)=\frac{1}{(2\pi)^{n/2}}|x|^nH_{2,0}^{1,2}\left[\frac{|x|^2}{2t^{1/2}}\left|
\begin{array}{c}
$$(1+(1-\alpha)j, 1/2)$$ \\
$$(n/2, 1/2), (2+j, 1/2) $$ \\
\end{array}\right.\right] ,\eqno{(3.12)}$$
$$h_{1j}(x,t)=\frac{1}{(2\pi)^{n/2}}|x|^nH_{2,0}^{1,2}\left[\frac{|x|^2}{2t^{1/2}}\left|
\begin{array}{c}
$$(1+(1-\alpha)j, 1/2)$$ \\
$$(n/2, 1/2), (1+j, 1/2) $$ \\
\end{array}\right.\right]. \eqno{(3.13)}$$
Computing the inverse transform of equation (3.9) and substituting into resulting equality formulas (3.12), (3.13), we get
$$G(x, t)=\frac{1}{(2\pi)^{n/2}|x|^n}\sum_{j=0}^{\infty}\frac{\left(-t^{1-\alpha}\right)^j}{j!}\left\{H_{2,0}^{1,2}\left[\frac{|x|^2}{2t^{1/2}}\left|
\begin{array}{c}
$$(1+(1-\alpha)j, 1/2)$$ \\
$$(n/2, 1/2), (2+j, 1/2) $$ \\
\end{array}\right.\right]+\right.$$
$$\left.+H_{2,0}^{1,2}\left[\frac{|x|^2}{2t^{1/2}}\left|
\begin{array}{c}
$$(1+(1-\alpha)j, 1/2)$$ \\
$$(n/2, 1/2), (1+j, 1/2) $$ \\
\end{array}\right.\right]\right\}.$$

Continuing to convert the equality (3.7) we can write formally
$$L^{-1}\left[\frac{1+s^{\alpha-1}}{s+s^{\alpha-1}+|\xi|^2}\tilde{g}(\xi)\right]=L^{-1}\left[F\left[G(x, s)\right](\xi)\tilde{g}(\xi)\right]. \eqno{(3.14)}$$
In view of (3.6) and (3.14), applying an inverse Laplace transform to equation (3.2) we finally obtain
$$\tilde{u}(\xi, t)=\sum_{j=0}^{\infty}(-1)^j\left(\mathcal{E}_{1, (1-\alpha)j+1,-|\xi|^2; 0+}^{j+1}\tilde{f}\right)(\xi,\tau)+F\left[G(x, t)\right](\xi)g(\xi) . \eqno{(3.15)}$$
To equation  (3.15) can  be further applied inverse Fourier transform and Fourier convolution property in sequence. Accordingly, the Theorem 1 is proven.

\section{The integro-differential diffusion equation with the Mittag-Leffler function in the kernel}

In this section we show equivalence of one  integro-differential diffusion equation with
the Mittag-Leffler function in the kernel to the  anomalous diffusion equation.

{\bf Theorem 2.} The integro-differential  diffusion equation

$$u_t-\triangle u+\int\limits_0^tk(t-\tau)\triangle u(x, \tau)d\tau=0, \ x\in\mathbb{R}^n, \ t>0 \eqno(4.1)$$
with memory
$$k(t)=t^{-\alpha}E_{1-\alpha, \ 1-\alpha}\left(-t^{1-\alpha}\right), \ \alpha\in(0, 1), $$
is equivalent to the time-fractional diffusion equation
$$u_t+\,_0^CD_t^{\alpha}u-\triangle u(x,t)=0. \eqno(4.2)$$

{\bf Proof.} Considering equation (4.1) as the Volterra integral
equation of the second kind with respect to $\triangle u$ for
fixed $x$ and applying Lemma 4, we have
$$\triangle u=u_t+\int\limits_0^tr(t-\tau) u_{\tau}(x, \tau)d\tau,  \eqno(4.3)$$
where $r(t)$ is resolvent of $k(t)$ and it satisfies the integral equation
$$r(t)=k(t)+\int\limits_0^tk(t-\tau)r(\tau)d\tau. \eqno(4.4)$$
In \cite{42} it was shown the resolvent of $t^{-\alpha}/\Gamma(1-\alpha)$ is the function $(d/dt)E_{1-\alpha}\left(-t^{1-\alpha}\right).$
Computing the derivative of Mittag-Leffler function \cite{43}, we get
$$\frac{d}{dt}E_{1-\alpha}\left(-t^{1-\alpha}\right)=t^{-\alpha}E_{1-\alpha, 1-\alpha}\left(-t^{1-\alpha}\right).$$
We apply to both sides of (4.4) the Laplace and denoting by $K(s)$ and $R(s)$ the imagines of origins $k(t)$ and $r(t),$ respectively, obtain
$$R(s)=K(s)+K(s)R(s). \eqno(4.5)$$
In view of $K(s)=L\left[t^{-\alpha}E_{1-\alpha, 1-\alpha}\left(-t^{1-\alpha}\right)\right]=1/\left(s^{1-\alpha}+1\right), \ \mathfrak{R}(s)>1$
(see \cite{44}), from  equality (4.5) we get
$$R(s)=\frac{K(s)}{1-K(s)}=\frac{1}{s^{1-\alpha}.}$$
From here it follows
$$L\left[R(s)\right]=L^{-1}\left[\frac{1}{s^{1-\alpha}}\right]=\frac{t^{-\alpha}}{\Gamma(1-\alpha)}=r(t).
$$
Then, (4.3) yields (4.2).

{\bf Remark.} Thus, the equation (4.1) with memory kernel $k(t)=t^{-\alpha}E_{1-\alpha, \ 1-\alpha}\left(-t^{1-\alpha}\right)$ describes the anomalously
diffusive transport of solute in heterogeneous porous media.

From this remark it follows that the solution of equation (4.1) with  conditions (1.2) can be given by formula (3.1) for $f(x, t)=0.$

\section{ Conclusions}
In applications, using the different types of a
memory kernel $k(t)$ in equation (4.1) it can be described a wide
variety of physical phenomena with memory effects. In this work,
it is shown that $n-$dimensional anomalous diffusion equation
(2.1) with $f(x,  t) = 0$  can be derived from the parabolic
integro-differential equation (4.1) with memory kernel
$t^{-\alpha}E_{1-\alpha, 1-\alpha}(-t^{1-\alpha}), \ \alpha\in(0,
1).$ Based on the Laplace transform method to the time variable
and Fourier transform to the spatial variable, the explicit
solution of initial–boundary problem for anomalous diffusion
equation is obtained. This solution includes the Prabhakar
fractional integral and Fox H-functions.

\end{document}